\documentclass[amstex]{article}
\usepackage{graphicx}
\usepackage{dcolumn}
\usepackage{amsmath}
\usepackage{amssymb}
\usepackage{amsfonts}
\usepackage{amscd}
\usepackage{bbm} %% for bb 1
\usepackage[usenames]{color}
\newtheorem{theorem}{Theorem}
\newtheorem{lmm}[theorem]{Lemma}
\newtheorem{cor}[theorem]{Corollary}
\newtheorem{pro}[theorem]{Proposition}
\newtheorem{df}[theorem]{Definition}
\newtheorem{rmk}[theorem]{Remark}

\newcommand\bfZ{{\bf Z}}
\newcommand\bfC{{\bf C}}
\newcommand\bfN{{\bf N}}
\newcommand\bfR{{\bf R}}

\newcommand\dprime{{\prime\prime}}
\newcommand\evl{\alpha_t}

\begin{document}
%%%%%%%%%%%%%%%%%%%%%%%%%%%%%%%%%%%%%%%%%%%
%%%%%%%%%%%%%%%%%%%%%%%%%%%%%%%%%%%%%%%%%%%
\newpage\thispagestyle{empty}
\begin{center}
{\huge\bf
Split Property 
\\
and 
\\
Fermionic String Order. }
\\
\bigskip\bigskip
\bigskip\bigskip
{\Large Taku Matsui}
\\
 Graduate School of Mathematics, Kyushu University,
\\
744 Motoka, Nishi-ku, Fukuoka 819-0395, JAPAN
\\
 matsui@math.kyushu-u.ac.jp
\\
\end{center}
\bigskip\bigskip\bigskip\bigskip
\bigskip\bigskip\bigskip\bigskip
{\bf Abstract:}  
We show split property of  gapped ground states for Fermion systems on a one-dimensional lattice and clarify
mathematical meaning of string oder of fermions.

%This is also for double spacing
\newpage
%%%%%%%%%%%%%%%%%%%%%%%%%%%%%%%%%%%%%%%%%%%
%%%%%%%%%%%%%%%%%%%%%%%%%%%%%%%%%%%%%%%%%%%
\section{Introduction}\label{Preliminary}
\setcounter{theorem}{0}
\setcounter{equation}{0}
Analysis of gapped ground states is a subject of intensive study in mathematical physics. Matrix product states
are typical examples for investigation of various aspects of gapped ground states. When we proceed to
study of more general gapped ground states for one-dimensional systems, we arrive at states of quantum spin chains 
with split property.Split property is statistical independence of two subsystems in an infinite  quantum spin chain. 
Historically the notion of split property was introduced in research of  local quantum field theories in the last century 
and it is relatively recent to apply it to analysis of ground states of quantum spin chains.
 In \cite{Matsui2000} and in  \cite{Matsui2010} we have shown that split property is valid for  gapped ground states
 of  infinite quantum spin chains and if split property is valid,
 the state has a matrix product state representation with an infinite dimensional auxiliary  boundary space.
 As the matrix product representation is constructed canonically on the basis of  the physical Hilbert space of the system,
 symmetry property of a gapped ground state is inherited to its  auxiliary  boundary space.
 Using matrix product representations, we have a Lieb-Mattis-Schultz type theorem for general SU(2) symmetric ground states. 
 Our result is extended to certain discrete symmetry including space reflection by H.Tasaki and Y.Ogata \cite{TasakiOgata2019}. 
 Split property is a standing base of $\bfZ_{2}$ invariant of Y.Ogata in \cite{Ogata2019}) which is a generalization of various invariant
 of symmetry protected topological phases.  
 \par
 The aim of the present paper is to investigate some basic aspects of split property for Fermion systems and
to clarify mathematical meaning of string order for quasi-one dimensional Fermion systems .
 \par
Though a part of results stated below can be proved in a more general setting of  unital $Z_{2}$ graded infinite dimensional $C^{*}$-algebras,
 for concreteness, we focus on CAR algebras( algebras generated by creation and annihilation operators of fermions).
 \par
 Let $\mathfrak A (K)$  be the (self-dual) CAR(canonical anticommutation relations) algebra $\mathfrak A (K) $ over a one-particle Hilbert space $K$ equipped with a conjugate unitary involution  $J$ . (c.f. \cite{Araki1971}.)
 Hereafter we assume that $K$ is separable.
 By the conjugate unitary involution, we mean $J$ is a conjugate linear operator acting on the complex Hilbert space $K$ satisfying
 \begin{equation}
 \label{eqn:a1}
J^{2} =1  , \quad ( J h , J k )_{K} = (k , h )_{K}   \quad ( h,k \in K) 
\end{equation}
 where $(k , h )_{K}$ is the inner product of vectors $k, h$ in $K$. 
 The CAR algebra $\mathfrak A (K) $ is the UHF $C^{*}$ algebra generated by $B(f) , f \in K$ where $B(f)$ depends linearly on $f \in K$ and satisfies the following canonical anticommutation relation and ${}^{*}-$relations:
 \begin{equation}
 \label{eqn:a2}
 B(h)^{*} = B(Jh) , \quad \{ B(h) , B(k)^{*} \} =  (k, h )_{K} 1.
 \end{equation}
  $\mathfrak A (K)$ is $Z_{2}$ graded where the $Z_{2}$ grading is specified with a parity automorphism $\Theta$  
determined via the following equations:  $\Theta (B(f)) = - B(f)$. Obviously, $\Theta^{2}(Q) = (Q)$ for any Q in $\mathfrak A (K) $.
The even (odd ) part of   $\mathfrak A (K)$ is denoted by $\mathfrak A^{(+)}(K)$ (resp.  $\mathfrak A^{(-)}(K)$)
\begin{equation}
 \label{eqn:a3}  
 \mathfrak A^{(\pm )}(K) = \left\{  \:\: Q \in  \mathfrak A (K) \:\: | \:\: \Theta (Q) = \pm Q   \:\: \right\}.
 \end{equation}
We presume that the one-particle Hilbert space $K$ is a direct sum of two $J$ invariant closed subspaces $K_{R}$ and $K_{L}$,
\begin{equation}
 \label{eqn:a4} 
 K = K_{R} \oplus K_{L}, \quad K_{R} \perp K_{L} , \quad JK_{R}=K_{R},\quad JK_{L}=K_{L}
 \end{equation}
and  to simplify notations, we set
\begin{eqnarray}
 \label{eqn:a5}
&& \mathfrak A = \mathfrak A (K), \quad
\mathfrak A_{L} = \mathfrak A (K_{L}), \:\:  \mathfrak A_{R} = \mathfrak A (K_{R}) ,
\nonumber
\\
&&\mathfrak A^{(\pm )}_{L} =   \mathfrak A^{(\pm )}(K) \cap \mathfrak A_{L} ,\:\:
\mathfrak A^{(\pm )}_{R} =   \mathfrak A^{(\pm )}(K) \cap \mathfrak A_{R} 
\end{eqnarray}
where $\mathfrak A (K_{R})$(resp. $\mathfrak A (K_{L})$) is the $C^{*}$ subalgebra of $\mathfrak A (K) $ generated 
by $B(f)$ with $f \in K_{R}$  (resp.$f \in K_{L}$). 
Note $\{ Q , R \} = 0$ for $Q \in\mathfrak A^{(-)}_{L} , R\in\mathfrak A^{(-)}_{R} $.
\par
For later use, we introduce another automorphsim $\Theta_{-}$ of $\mathfrak A$ 
 determined via the following equations:
\begin{eqnarray}
 \label{eqn:a6}
 &&\Theta_{-}(B(f)) = - B(f) , \:\:\Theta_{-}(B(g)) = B(g) \quad ( f \in K_{L},  \:\: g \in K_{R}).
\end{eqnarray}
Let $\hat{\mathfrak A}$ be the crossed product of the $\bfZ_{2}$ action induced by $\Theta_{-}$, 
$\hat{\mathfrak A} =  \mathfrak A \rtimes \bfZ_{2}$, namely
$\hat{\mathfrak A}$ is generated by $\mathfrak A$ and a self-adjoint unitary $T$ implementing
$\Theta_{-}$ 
\begin{equation}
 \label{eqn:a601}
  T^{2}=1 , \quad T=T^{*},\quad TQT^{-1} = \Theta_{-}(Q) , \quad Q \in \mathfrak A
\end{equation}  
We extend $\Theta$  via the equation $\Theta (T) = T$  to  an automorphism  of  $\hat{\mathfrak A}$.
For any $Q \in  \hat{\mathfrak A}$ we set 
\begin{equation}
 \label{eqn:a602}
  Q_{\pm} = 1/2 ( Q \pm \Theta (Q) ) . 
\end{equation}  
Let   $\mathfrak A^{S}$ be the $C^{\ast}$-subalgebra of $\hat{\mathfrak A}$ generated by
$\mathfrak A_{L}$ and  $T\mathfrak A_{R}^{(-)} $ and let
$\mathfrak A^{S}_{R}$ be the $C^{\ast}$-subalgebra of $\mathfrak A^{S}$ generated by
$T\mathfrak A_{R}^{(-)} $ and we set 
$$\mathfrak A^{S \; (\pm )} = \{ Q \in \mathfrak A^{S} \: | \: \Theta (Q) = \pm Q \} ,\quad
\mathfrak A^{S \; (\pm )}_{L/R} = \mathfrak A^{S \; (\pm )} \cap \mathfrak A^{S}_{L/R}. $$
Obviously,  $\mathfrak A^{S \; (+)}_{L/R} = \mathfrak A^{(+)}_{L/R}$. 
$\mathfrak A_{L}$ and  $T\mathfrak A_{R}^{(-)} $ commute and
$\mathfrak A_{L/R}$,  and $ \mathfrak A^{S}_{L/R}$ are isomorphic to the UHF $C^{*}$-algebra of  the type $2^{\infty}$,
 as a consequence, we can identify  $\mathfrak A^{S}$  with $\mathfrak A_{L} \otimes \mathfrak A_{R}$:
\begin{equation}
 \label{eqn:a7}
 \mathfrak A^{S} =  \mathfrak A_{L} \otimes \mathfrak A_{R} \quad \mathfrak A^{S}_{L}=  \mathfrak A_{L} \otimes 1
,\quad   \mathfrak A^{S}_{R}=1 \otimes  \mathfrak A_{R}.
\end{equation}
Let us recall definition of  ($\bfZ_{2}$ graded) product states.  Note that our definition of $\mathfrak A^{S}$ 
is different from that of the Pauli spin system via Jordan Wigner transformation in  \cite{Araki1971}. 
The definition is convenient for proof of split property in Section \ref{Split}.
\begin{lmm}
\label{lmm:1.1}
% definition of  $Z_{2}$ graded product states
Let $\psi_{L}$ (resp. $\psi_{R}$) be a state of  $\mathfrak A_{L}$ (resp. $\mathfrak A_{R}$) and suppose that
$\psi_{L}$  is $\Theta$ invariant.
There exists a state $\psi_{L}\otimes_{\bfZ_{2}} \psi_{R}$ of $\mathfrak A$ satisfying
\begin{equation}
\label{eqn:a8}
\psi_{L}\otimes_{\bfZ_{2}}\psi_{R}(QP) = \psi_{L}(Q_{+}) \psi_{R}(P)
\end{equation}
fro any $Q \in \mathfrak A_{L}$ and any $P \in \mathfrak A_{R}$.
\end{lmm}
  For proof of this lemma, it suffices to present the GNS representation associated with 
$\psi_{L}\otimes_{\bfZ_{2}} \psi_{R}$. Let $\{ \pi_{L}(\mathfrak A_{L} ), \Omega_{L}, \mathfrak H_{L} \}$
be the GNS triple of $\mathfrak A_{L}$ associated with the state $\psi_{L}$
where  $\pi_{L}(\cdot )$ is the representation of $\mathfrak A_{L}$ on the Hilbert space $\mathfrak H_{L} $ with
the GNS cyclic vector $\Omega_{L}$, and let  $\{ \pi_{R}(\cdot ), \Omega_{R}, \mathfrak H_{R} \}$
be the GNS triple of $\mathfrak A_{R}$ associated with the state $\psi_{R}$.
Since $\psi_{L}$  is $\Theta$ invariant, there exists a self-adjoint unitary $\Gamma_{L}$ acting on  $\mathfrak H_{L} $
satisfying
\begin{equation}
\label{eqn:a801}  
\Gamma_{L} \Omega_{L} = \Omega_{L} , \:\: \Gamma_{L}  \pi_{L}(Q)\Gamma_{L}^{-1} = \pi_{L}(\Theta (Q)) ,\quad  Q \in \mathfrak A_{L}  .
\end{equation}
Let $P_{L}$ be the orthogonal projection from $K$ to $K_{L} $ and  $P_{R}$ be that to $K_{L} $ and set
\begin{equation}
\label{eqn:a802}
\pi (B(h))  = \pi_{L} (B(P_{L} h) \otimes 1 + \Gamma_{L} \otimes \pi_{R} (B( P_{R}h)  \:  h \in K .
\end{equation}
Then, (\ref{eqn:a802}) gives rise to the GNS representation of $\mathfrak A$ associated with  $\psi_{L}\otimes\psi_{R}$.
\bigskip
\newline
If , in stead, $\psi_{R}$  is $\Theta$ invariant, we have a state  $\psi_{L}\otimes_{\bfZ_{2}}\psi_{R}$ determined by 
\begin{equation}
\label{eqn:a9}
\psi_{L}\otimes_{\bfZ_{2}}\psi_{R} (QP) = \psi_{L}(Q) \psi_{R}(P_{+} )
\end{equation}
fro any $Q \in \mathfrak A_{L}$ and any $P \in \mathfrak A_{R}$.
We employ the same notation  $\psi_{L}\otimes_{\bfZ_{2}} \psi_{R}$ for the state
determined by (\ref{eqn:a8}) or  by (\ref{eqn:a9}) and  $\psi_{L}\otimes_{\bfZ_{2}} \psi_{R}$ will be referred to as the graded product state of 
$\psi_{L}$ and $\psi_{R}$ or simply the product state if there is no risk of confusion.
% definition of  split property
\begin{df}
\label{df:1.2}
Let $\psi$ be a state of $\mathfrak A$. We say split property is valid for $\psi$ if
$\psi$ is quasi-equivalent to a product state $\psi_{L}\otimes_{\bfZ_{2}} \psi_{R}$.
\end{df}
In the above Definition \ref{df:1.2}, we are assuming, at least, one of  $\psi_{L}$ and $\psi_{R}$ is $\Theta$ invariant. 
If a state $\psi$ of  $\mathfrak A$ is $\Theta$ invariant, satisfying split property, we may suppose that 
both $\psi_{L}$ and $\psi_{R}$ are $\Theta$ invariant, for example, if $\psi_{L}$  is $\Theta$ invariant and
$\psi$ is quasi-equivalent to a product state $\psi_{L}\\otimes_{\bfZ_{2}} \psi_{R}\circ\Theta$,
hence,  $\psi$ is quasi-equivalent to  $\psi_{L}\otimes_{\bfZ_{2}} 1/2(\psi_{R}+ \psi_{R}\circ\Theta )$ as well.
\begin{theorem}
\label{th:1.3}
Let $\psi$ be a pure state of  $\mathfrak A$ which is $\Theta$ invariant. 
\\
Suppose that split property is valid for $\psi$.
Then, $\psi$ restricted to $\mathfrak A_{R}$ (resp. $\mathfrak A_{L}$) gives rise to a type I representation of $\mathfrak A_{R}$.
(resp. of $\mathfrak A_{L}$).
\\
Conversely, if $\mathfrak A_{R}$ (resp. $\mathfrak A_{L}$) gives rise to a type I representation, the split property is valid
for  $\psi$.
\end{theorem}
Next we turn to gapped ground states.
We will see that gapped ground states with string order
cannot be connected to the Fock vacuum state in Section \ref{String}.  
\par
Let $\mathfrak A$ be the CAR algebra generated by
creation and annihilation operators $c^{*}_{j} , c_{i}\:\: i,j \in\bfZ$ satisfying the standard canonical anticommutation relations:
\begin{equation}
\label{eqn:a10}
\{ c_{j} , c_{i} \} =0 ,\:\:\: \{  c^{*}_{j} , c^{*}_{i} \} = 0,\:\:\  \{ c^{*}_{j} , c_{i} \} = \delta_{i,j} 1 . 
\end{equation}
The relation to the selfdual formalism of . \cite{Araki1971}  will be explained later in Section \ref{String}.
%$K = l^{2}(\bfZ )\oplus l^{2}(\bfZ )$
\par
For integers $n < m$, let  $\mathfrak A_{[n,m]}$ be the subalgebra of $\mathfrak A$ generated by  
$c^{*}_{j} , c_{i}$   ( $n \leq i,j \leq m $)
and $\mathfrak A_{[n,m]}^{(+)}$ be the even part of $\mathfrak A_{[-n,n]}$. 
Let $\mathfrak A_{loc}$ be the algebra of strictly local observables:
\begin{equation}
\label{eqn:a11}
 \mathfrak A_{loc}  \:\:\: = \:\:\: \cup_{-\infty < n < m <\infty} \mathfrak A_{[n,m]}. 
\end{equation}
By $\tau_{k}$ ($k \in \bfZ$) we denote the lattice translation
which is an automorphism of $\mathfrak A (K)$ satisfying 
\begin{equation}
\label{eqn:a12}
\tau_{k}(c_{i} ) = c_{i+k}, \tau_{k}(c_{i}^{*} ) = c_{i+k}^{*}.
\end{equation}
Fix $h = h^{*} \in  \mathfrak A_{[-n,n]}^{(+)}$ for some $n >1$ and we consider a finite volume Hamiltonian $H_{M}$
defined by
\begin{equation}
\label{eqn:a13}
  H_{M} = \sum_{-M+n \leq k \leq M-n}   \tau_{k}(h) \quad \in \mathfrak A_{[-M,M]}^{(+)} .
\end{equation}
The formal infinite volume Hamiltonian 
\begin{equation}
\label{eqn:a13-2}
H = \lim_{M \to\infty} H_{M}
\end{equation}
 generates a time evolution of  our Fermion system,
more precisely,  an infinite volume time evolution  $\evl^{h} (Q)$ of an observable $Q$ is determined via the following equation:
\begin{equation}
\label{eqn:a14}
 \lim_{M\to\infty} \:\: e^{it H_{M} } Q  e^{-it H_{M}}  =  \evl^{h} (Q) \quad Q \in \mathfrak A (K) .
\end{equation}
\begin{df}
\label{bf:ground state}
Let $\psi$ be a translationally invariant state. 
\\
(i) $\psi$ is a ground state of $H$  if $\psi$ minimizes the energy density:
\begin{equation}
\label{eqn:a15}
   \psi (h )  \:\:\: = \:\:\: \inf \varphi (h)
\end{equation}
where $\inf$ is taken among the set of all translationally invariant states.
\\
(ii) A translationally invariant ground state $\psi$ is gapped if the following inequality is valid for some positive constant $m$:
\begin{equation}
\label{eqn:a16}
\psi (Q^{*} [H,Q]) \geq m (\psi (Q^{*}Q) -  |\psi (Q)|^{2}) \quad Q \in  \mathfrak A_{loc} .
\end{equation}
\end{df}
\begin{rmk}
\label{rmk:1.6}
(i) Any translationally invariant pure state $\psi$ is $\Theta$ invariant. This is a folklore among specialists, and we present a brief
 sketch of proof. Suppose  $\psi$ is pure, translationally invariant but not $\Theta$ invariant. 
 There exists $Q \in \mathfrak A^{(-)}$ and an increasing sequence of
even numbers, $j_{k} \in 2 \bfN$ such that
$j_{1}< j_{2} < j_{k} < j_{k+1} \dots$ such that
$$\psi (Q) \ne 0 , \quad \lim_{k\to\infty} \pi_{\psi}(\tau_{j_{k}}(Q)) = W \ne 0 .$$
As a consequence, we have
$$ W \pi_{\psi}(Q)=\pi_{\psi}(\Theta (Q))W, \:\: W^{*} \pi_{\psi}(Q)=\pi_{\psi}(\Theta (Q))W^{*} , \:\:W^{2} =c 1, $$
for some $\:\: c \in \bfC$.
By multiplying a suitable constant, we may assume that $W$ is a self-adjoint unitary $W=W^{*},\: W^{2}=1$ implementing $\Theta$.
If we denote the normal extension of $\tau_{1}$ to $\pi_{\psi}(\mathfrak A)^{\dprime}$ by the same symbol, 
we have $\tau_{1}(W)= \pm W$. This is because $\tau_{1}(W)$ implements $\Theta$ as well.
Set $\xi_{\pm} =\frac{1}{2}(1 \pm W)\Omega_{\psi}$. The positive normal functional $\psi_{\pm}(Q)=(\xi_{\pm},\pi_{\psi}(Q)\xi_{\pm})$
are $\Theta$ and $\tau_{2}$ invariant because
$$\psi_{\pm}(\tau_{2}(Q)) = \frac{1}{4} ((1 \pm \tau_{-2}(W))\Omega_{\psi}, \pi_{\psi}(Q) (1 \pm \tau_{-2}(W))\Omega_{\psi})$$
$$\psi_{\pm}(\Theta (Q)) = \frac{1}{4} (W (1 \pm W)\Omega_{\psi}, \pi_{\psi}(Q) W (1 \pm W)\Omega_{\psi} ).$$
By definition, $\psi_{\pm}(W) =\pm1$, however, due to $\Theta$ and $\tau_{2}$ invariance of $\psi_{\pm}$
$$ \psi_{\pm}(\tau_{j_{k}}(Q)) = 0$$ 
because $Q \in \mathfrak A^{(-)}$. $W=\lim_{k\to\infty}\pi_{\psi}(\tau_{j_{k}}(Q))$ implies a contradiction. 
\end{rmk}
\begin{rmk}
\label{rmk:1.7}
If a translationally invariant state $\psi$ satisfies (\ref{eqn:a15}), $\psi$ is $\evl^{h}$ invariant. 
On the GNS Hilbert space $\mathfrak H_{\psi}$ with the GNS cyclic vector $\Omega_{\psi}$, there exists a positive self-adjoint operator
$H_{\psi}$ on $\mathfrak H_{\psi}$ satisfying 
$$  e^{it H_{\psi} } \pi_{\psi}(Q)\Omega_{\psi} = \pi_{\psi}(\evl (Q))\Omega_{\psi} .$$
$H_{\psi}$ is regarded as the regularized Hamiltonian.
\par
By $spec \: H_{\psi}$ iwe denote the spectrum of $H_{\psi}$.
The condition (\ref{eqn:a16}) is equivalent to the spectrum gap condition, namely, the infimum of  the spectrum of $H_{\psi}$
is non-degenerate and
$spec \: H_{\psi} \subset \{ 0 \} \cup [ m , \infty )$.
\end{rmk}
\begin{df}
\label{bf:1.5}
Let $\psi$ be a translationally invariant pure state.  $\psi$ has string order if there exist $Q_{1}$in 
$\mathfrak A^{(-)} \cap \mathfrak A_{[n , -1]}$ and 
$Q_{2}$ in $\mathfrak A^{(-)}_{loc} \cap \mathfrak A_{[0, m]}$ with $n < 0 \leq m$ 
satisfying
\begin{equation}
\label{eqn:a17}
\lim_{k \to \infty} \psi ( Q_{1} \: S[0,2k-1]\: Q_{2} ) \ne 0
\end{equation}
where
$$ S[0,k-1]  = \prod_{j=0}^{k-1}\:\:\: (2 c^{*}_{j}c_{j} - 1).$$
\end{df}
By the standard Fock state $\psi_{F}$ we mean  a state of $\mathfrak A (K)$ uniquely 
determined by the equation  $\psi_{F}(c^{*}_{j} c_{j}) = 0$ for any $j \in \bfZ$.
\begin{theorem}
\label{th:1.5}
Let $\psi$ be a translationally invariant pure gapped ground state for $H$.
If $\psi$ has string order, $\psi$  cannot be connected to the standard Fock state $\psi_{F}$
by any one-parameter group of automorphisms $\evl^{h^{\prime}}$,
\begin{equation}
\label{eqn:a18}
 \psi \circ \evl^{h^{\prime}} \ne \psi_{F} 
\end{equation} 
for any $t$ and any $h^{\prime} \in \mathfrak A_{loc}$.
\end{theorem}
\begin{rmk}
\label{rmk:1.9}
Historically, in \cite{denNijs1989},  M.den Nijs, and K.Rommelse introduced 
the string order of gapped ground states for quantum spin chains 
 to characterize the Haldane phase.
\par
In their pioneer work \cite{Kennedy1992} of the Haldane phase, T. Kennedy and H. Tasaki discovered a relation between
the string order of  M.den Nijs, and K.Rommelse and hidden $\bfZ_{2}\times\bfZ_{2}$ symmetry breaking.
We recommend a monograph  \cite{Tasaki2020} by H.Tasaki for  detail and further extension of research.  
\par
 Turning to string order of Fermion systems discussed here, among various physics literature, and in the present context,
 one relevant is  a paper  \cite{Y.Bahri2014} of Y. Bahri and A. Vishwanath 
 Y. Bahri and A. Vishwanath examined  Fermion systems associated with the XY model. 
 According to Y. Bahri and A. Vishwanath, string order plays a role of detecting majorana fermions.
The string order we consider here is same as that of Y. Bahri and A. Vishwanath and 
it characterizes the hidden $\bfZ_{2}$ symmetry breaking 
where the hidden system is a Pauli spin chain obtained by Jordan Wigner transformation. 
Our results are valid for general short range  periodic Hamiltonians with gapped ground states, in particular, 
quasi one-dimensional tight binding models are within reach of our results.
\end{rmk}
\begin{rmk}
\label{rmk:1.10}
Theorem \ref{th:1.5} is valid for periodic gapped ground states. Theorem \ref{th:1.5} is valid
not only for finite range Hamiltonians but for long range Hamiltonians as long as the interaction across
$(-\infty , -1]$ and $[0 \infty )$ is bounded.
\par
We expect the stability of the Fermion string order can be proven for states change under 
one-parameter family of automorphisms considered in \cite{MoonOgata2020}, however, in this paper we concentrate on
general mathematical parts.
\par
In Theorem \ref{th:1.5} we consider translationally invariant states. We believe a similar result is valid for non-translationally invariant
setting. In that case, the bottom of spectrum can be infinitely degenerate, and some additional work is needed.
\end{rmk}
\bigskip
After completing our first manuscript, we noticed a work \cite{C.Bourne2019} of C.Bourne and
H. Schulz-Baldes . There are overlaps  both in methods and in contents between C.Bourne and H. Schulz-Baldes
ours. We mention here difference between the two.
\begin{description}
\item[(i)]  Theorem \ref{th:1.3} was remarked without proof in our previous paper \cite{Matsui2010}.
Consider the translationally invariant Hamiltonian for $\mathfrak A$ defined by 
\begin{equation}
\label{eqn:a19}
H^{F} = - \sum_{k= -\infty}^{\infty} \:\: (c_{k}^{*} - c_{k} )(c_{k+1}^{*} + c_{k+1} ) .
\end{equation}
After the Jordan Wigner transformation, the corresponding Hamiltonian $H^{P}$ for our Pauli spin algebra $\mathfrak A^{P}$ is
\begin{equation}
\label{eqn:a20}
H^{P} = - \sum_{k= -\infty}^{\infty} \:\: \sigma_{x}^{(k)}  \sigma_{x}^{(k+1)} .
\end{equation}
A unique translationally invariant , $\Theta$ invariant ground state $\psi^{P}$ of $H^{P}$ is an average of two product states 
$\varphi$ and $\varphi \circ\Theta$:
\begin{equation}
\label{eqn:a21}
\psi^{P} = \frac{\varphi + \varphi \circ\Theta}{2} , \quad \varphi ( \sigma_{x}^{(k)}) = 1  \quad k \in \bfZ .
\end{equation}
(\ref{eqn:a21}) shows that restriction of $\psi$ for $H^{F}$ to  $\mathfrak A_{L}^{(+)} \otimes \mathfrak A_{R}^{(+)}$ is not a factor state, not a product state.We have to take into account this possibility  in our proof of  Theorem \ref{th:1.3}.
\item[(ii)] We proved that any translationally invariant pure gapped ground states satisfy split property and the $\bfZ_{2}$ index
is well-defined for  those states. To prove our claim we employ results of   \cite{ArakiEvans1983} and \cite{ArakiMatsui1985}.
\item[(iii)] The definition of $\bfZ_{2}$ index introduced in \cite{C.Bourne2019} by C.Bourne and H. Schulz-Baldes and that by ourself
are different  but equivalent. We claim that the $\bfZ_{2}$ index is determined by presence or absence of string order.
\end{description}
We believe both  \cite{C.Bourne2019} of C.Bourne and H. Schulz-Baldes  and ours are complementary
and both papers shed a new scope of string order of Fermion chains and the $\bfZ_{2}$ index.
\par
We will prove results within the framework of the theory of operator algebras. Basic references of this approach
are  \cite{BratteliRobinsonI} and \cite{BratteliRobinsonII}. 
For more physics oriented readers, the monograph \cite{Tasaki2020} of H.Tasaki is a very good introduction 
to various mathematical aspects of quantum many body problems.
%%%%%%%%%%    End of Section 1        %%%%%%%%%%%%%%%%%%%
%%%%%%%%%%    End of Section 1        %%%%%%%%%
\section{Proof of Theorem \ref{th:1.3} }\label{Split}
\setcounter{theorem}{0}
\setcounter{equation}{0}
For later use, we collect results of $\Theta$  invariance of pure states.
Let $\mathfrak B$ be a unital $C^{*}$-algebra and  $\Theta$ be an involutive automorphism of $\mathfrak B$, 
$\Theta^{2}(Q) = Q$ for $Q \in \mathfrak B$. 
\\
Set 
$$\mathfrak B^{(\pm)} = \{ Q \in \mathfrak B \:\: | \:\: \Theta (Q) = \pm Q \:\:\} .$$
We assume that $\mathfrak B^{(-)}$ is generated by a single self-adjoint unitary $Z$ in $\mathfrak B^{(-)}$ and $\mathfrak B^{(+)}$,
$$ \mathfrak B^{(-)} = Z\mathfrak B^{(+)}$$
Let $\psi$ be a $\Theta$ invariant state of $\mathfrak B$. We denote the GNS triple of  $\mathfrak B$ associated with $\psi$ by
$\{ \pi_{\psi }(\mathfrak B) , \Omega_{\psi} , \mathfrak H_{\psi} \}$ where $\pi_{F}(\mathfrak B)$ is the representation of $\mathfrak B$ 
on the GNS Hilbert space $\mathfrak H_{\psi}$ and , $\Omega_{\psi} \in \mathfrak H_{\psi}$ is the GNS cyclic vector. 
As $\psi$ is $\Theta$ invariant, there exists a unique self-adjoint unitary operator $\Gamma (\Theta )$ on $\mathfrak H_{\psi}$ satisfying
$$\Gamma (\Theta )\pi_{\psi}(Q)\Omega_{\psi} = \pi_{\psi}(\Theta (Q)\Omega_{\psi} \:\:(Q \in \mathfrak B ),
\quad \Gamma (\Theta )^{*}=\Gamma (\Theta ),\:\:\Gamma (\Theta )^{2} = 1.$$
We set  $\mathfrak H_{\psi}^{\pm} =  \overline{ \pi_{\psi }(\mathfrak B^{(\pm )})\Omega_{\psi}}$ and the representation of
$\mathfrak B^{(+)}$ restricted to $\mathfrak H_{\psi}^{\pm}$ is denoted by $\pi_{\psi }^{(\pm)}(\cdot )$.
%%%%%  Lemma2.1  %%%%%%%%
\begin{lmm}[ Lemma 4.1 of \cite{ArakiMatsui1985}]
\label{lmm:2.1}
Let $\mathfrak B$ $\Theta$ be as above and let $\psi$ be a $\Theta$ invariant state of $\mathfrak B$. 
\\
(i)  Suppose that  $\psi$ is a pure state of $\mathfrak B$. Then, the restriction of $\psi$ to $\mathfrak B^{(+)}$ 
is a pure state of $\mathfrak B^{(+)}$.
\\
(ii)  Suppose that the  GNS representations  $\pi_{\psi }^{(+)}(\cdot )$ of $\mathfrak B^{(+)}$ is irreducible.
Then, $\psi$ is a pure state of $\mathfrak B$ if and only if 
$\pi_{\psi }^{(\pm)}(\cdot )$ of $\mathfrak B^{(+)}$ are mutually disjoint(= not unitarily equivalent) .
\\
(iii) Suppose that the  GNS representations  $\pi_{\psi }^{(+)}(\cdot )$ of $\mathfrak B^{(+)}$ is irreducible
and that   $\psi$ is not a pure state of $\mathfrak B$.
Then, there exists a pure state $\varphi$ of $\mathfrak B$ satisfying 
$\psi = \frac{ \varphi + \varphi\circ\Theta }{2}$.  $\varphi$  and $\varphi\circ\Theta$ are  mutually disjoint.
\end{lmm}
%%%%%%% Lemma2.2 %%%%%%%
\begin{lmm}
\label{lmm:30.9}
Let $\mathfrak B$ and $\Theta$  be as above.
Suppose that a state $\varphi$ of $\mathfrak B$ is pure and that $\varphi$ and $\varphi\circ\Theta$ are disjoint.
\\
Then, $\varphi$ restricted to $\mathfrak B^{(+)}$ is pure, and 
$\pi_{\varphi}(\mathfrak B)^{\dprime} = \pi_{\varphi}(\mathfrak B^{(+)})^{\dprime}$.
\end{lmm}
{\it Proof.}
Set 
$$\psi = \frac{1}{2} \left(\varphi + \varphi\circ\Theta \right) .$$
We employ the GNS triple $\{ \pi_{\psi}(\mathfrak B) , \Omega_{\psi} , \mathfrak H_{\psi} \}$
 associated with $\psi$ is realized in terms of representations associated with $\varphi$.
$$\mathfrak H_{\psi} =  \mathfrak H_{\varphi}\oplus  \mathfrak H_{\varphi} ,\quad 
\Omega_{\psi} = \frac{1}{\sqrt{2}} \left( \Omega_{\varphi}\oplus\Omega_{\varphi}\right), \quad
\pi_{\psi}(Q) = \pi_{\varphi}(Q) \oplus \pi_{\varphi}(\Theta (Q))$$
We now denote $h= (f,g)\in \mathfrak H_{\psi}$ by a column vector and any bounded operator on $\mathfrak H_{\psi}$ by
a 2 by 2 matrix as follows.
$$ h =  \begin{pmatrix} f \\ g \end{pmatrix} \in \mathfrak H_{\psi},\quad
\pi_{\psi}(Q) =  \begin{pmatrix} \pi_{\varphi}(Q) & 0\\ 0 &\pi_{\varphi}(\Theta (Q)) \end{pmatrix} , \quad Q \in \mathfrak B.
$$
Let $V(\Theta )$ and $X$ be self-adjoint unitaries on  $\mathfrak H_{\psi}$ defined via the following identities :
\begin{equation}
\label{eqn:c26}
V(\Theta ) \: = \: \left(\begin{array}{cc}0 & 1 \\1 & 0\end{array}\right) \in \mathfrak B(\mathfrak H_{\psi} ) , \quad
X = \left(\begin{array}{cc}1 & 0 \\0 & -1\end{array}\right).
\end{equation}
Due to disjointness of $\varphi$ and $\varphi\circ\Theta$, $X$ is an element of $\pi_{\psi}(\mathfrak B)^{\dprime}$ and
by a direct calculation, we see
\begin{equation}
\label{eqn:c27}
V(\Theta ) X  V^{-1}(\Theta ) = - X ,\quad V(\Theta ) \pi_{\psi}(Q) V^{-1} (\Theta ) =  \pi_{\psi}(\Theta (Q)), \quad Q \in \mathfrak B .
\end{equation}
The first identity in (\ref{eqn:c27}) suggests that $X$ is in the weak closure of $\pi_{\psi}(\mathfrak B^{(-)})$.
The center of the von Neumann algebra $\pi_{\psi}(B)^{\dprime}$ is two dimensional , generated by $X$ and
\begin{eqnarray*}
&&\varphi (Q) = \left( \Omega_{\psi} , \pi_{\psi}(Q) \frac{(1+X)}{2} \Omega_{\psi} \right) ,
\\
&&\varphi\circ\Theta (Q) = \left( \Omega_{\psi} , \pi_{\psi}(Q)\frac{(1-X)}{2}  \Omega_{\psi} \right),\quad Q \in B.
\end{eqnarray*}
We claim that 
\begin{equation}
\label{eqn:c28}
 \overline{ \pi_{\psi}(\mathfrak B^{(+)}) \Omega_{\varphi}}  = \mathfrak H_{\varphi} = \overline{\pi_{\psi}(\mathfrak B)\Omega_{\varphi}} .
\end{equation}
As the self-adjoint unitary $X$ is in in the weak closure of $\pi_{\psi}(\mathfrak B^{(-)})$, 
any element $R$ in  the weak closure of $\pi_{\psi}(\mathfrak B^{(-)})$
is a product 
\begin{equation}
\label{eqn:c29}
R = R^{+} X , \quad R^{+} = RX \in  \pi_{\psi}(\mathfrak B^{(+)})^{\dprime} .
\end{equation} 
We obtain
$$X \begin{pmatrix} \Omega_{\varphi}\\  0 \end{pmatrix} = \begin{pmatrix}\Omega_{\varphi} \\ 0 \end{pmatrix}  , \quad
R \begin{pmatrix} \Omega_{\varphi} \\ 0 \end{pmatrix} = R^{+} \begin{pmatrix} \Omega_{\varphi} \\ 0 \end{pmatrix}  $$
which implies (\ref{eqn:c28}).
\par
Next we show that $\varphi$ restricted to $\mathfrak B^{(+)}$ is pure. If this is not the case,
there exists a non-trivial projection $E$ in the commutant $\pi_{\varphi}(\mathfrak B^{(+)})^{\prime}$. 
As  $\Omega_{\varphi}$ is a cyclic vector for 
$\pi_{\varphi}(\mathfrak B^{(+)})^{\dprime}$,
$E\Omega_{\varphi} \ne 0$ in $\mathfrak H_{\varphi}$. Set
$$\tilde{E} = \begin{pmatrix} E & 0\\ 0 &E \end{pmatrix} \in \mathfrak B (\mathfrak H_{\psi}).$$
Then, $\tilde{E}$ commutes with $X$ , and with any element of $\pi_{\psi}(\mathfrak B^{(-)})$ due to (\ref{eqn:c29}). 
However, this contradicts with irreducibility of $\pi_{\varphi}(B)$ on $\mathfrak H_{\varphi}$. 
It turns out that $\varphi$ restricted to $\mathfrak B^{(+)}$ is pure.
As $\pi_{\varphi}(\mathfrak B^{(+)})^{\dprime}$ is trivial, 
$$\pi_{\varphi}(\mathfrak B^{(+)})^{\dprime} = \pi_{\varphi}(\mathfrak B)^{\dprime}=  \mathfrak B (\mathfrak H_{\varphi}) .$$
{\it End of Proof.}
%%%%%%
\bigskip
\\
We will borrow the following lemma in  \cite{Stratila1978}  due to S.Str$\breve{a}$til$\breve{a}$ and D.Voiculescu  as well.
(c.f. 2.7 of \cite{Stratila1978})
\begin{lmm}
\label{lmm:2.0}  
 Let $\mathfrak B$  and $\Theta$  be as in the previous lemma  $\ref{lmm:2.1}$ and $\psi_{1}$ and $\psi_{2}$
be $\Theta$ invariant states of $\mathfrak B$. If $\psi_{1}$ and $\psi_{2}$ restricted to $\mathfrak B^{(+)}$
are quasi-equivalent, so are $\psi_{1}$ and $\psi_{2}$.
\end{lmm}
%%%%%
\begin{lmm}
\label{lmm:2.2}
(i) Suppose that  $\varphi$ is a pure state of $\mathfrak A^{S}$.  
 Restriction of $\varphi$  to $\mathfrak A^{S}_{L}$ and that  to $\mathfrak A^{S}_{R}$ are factor states.
\\
(ii) Suppose that  $\psi$ is a $\Theta$ invariant factor state of $\mathfrak A$.  
\\
(ii-a) If $\psi$ is not quasi-equivalent to  $\psi\circ\Theta_{-}$, restriction of $\psi$ to $\mathfrak A_{L}$ and that to 
$\mathfrak A_{R}$ are factor states.
\\
(ii-b) If $\psi$ is pure , equivalent to  $\psi\circ\Theta_{-}$,  and the restriction of $\psi$ to $\mathfrak A^{(+)}$ and that of 
$\psi\circ\Theta_{-}$ to $\mathfrak A^{(+)}$ are equivalent, then, the restriction of $\psi$ to $\mathfrak A_{L}$ and $\psi$ restricted 
 to $\mathfrak A_{R}$ are factor states.
 \\
(ii-c)  If $\psi$ is equivalent to  $\psi\circ\Theta_{-}$,  and $\psi$ restricted to $\mathfrak A^{(+)}$ and $\psi\circ\Theta_{-}$
 restricted to $\mathfrak A^{(+)}$ are not equivalent, then, $\psi$ restricted to $\mathfrak A_{L}$ and
 to $\mathfrak A_{R}$ are not factor states. $\psi$ restricted to $\mathfrak A_{L}$ is an average of two mutually non quasi-equivalent
 factor states.
 \end{lmm}
 Examples of Lemma \ref{lmm:2.2} (ii) are provided by  Fock states associated with
 ground states of one-dimensional XY model. (c.f. \cite{ArakiMatsui1985})
 \\
 \\
{\it Proof of Lemma \ref{lmm:2.2}}
\\
The claim of Lemma \ref{lmm:2.2} (i) follows from identification of the center of $\pi_{\varphi}(\mathfrak A^{S})^{\dprime}$
and the algebra at infinity. (c.f. Theorem 2.6.10 of \cite{BratteliRobinsonI}.) 
The claim of Lemma \ref{lmm:2.2} (ii) is not mentioned in this form in \cite{BratteliRobinsonI} and we present our proof here.
\par
First we recall proof of   Lemma \ref{lmm:2.2} (i)  briefly.
As $C^{*}$-algebras, $\mathfrak A_{L}^{S}$ and $\mathfrak A_{R}^{S}$ are an infinite tensor product of $M_{2}(\bfC )$. 
We assign a non-negative integer to specify each tensor product component of $\mathfrak A_{R}^{S}$ 
and a negative  integer to specify each tensor product component of $\mathfrak A_{R}^{L}$. 
$\mathfrak A_{R}^{S}$ is generated by Pauli spin matrix $\sigma^{(j)}_{\alpha}$
where $j \in \bfZ , 0\leq j$ ,$\alpha =x,y,z$ and $\mathfrak A_{L}^{S}$ is generated by Pauli spin matrix $\sigma^{(j)}_{\alpha}$
where $j \in \bfZ ,  j < 0$ ,$\alpha =x,y,z$. 
\par
Let $\Lambda$ be a subset of $\bfZ$ and $\mathfrak A^{S} (\Lambda )$ be the
$C^{*}$-subalgebra of $\mathfrak A^{S}$ generated by $\sigma^{(j)}_{\alpha}$ $j \in \Lambda , \alpha =x,y,z$.
$\mathfrak A^{S} (\Lambda )$ is the set of observables localized in  $\Lambda$.
We set  
$$\mathfrak A^{S}_{loc} = \cup_{\Lambda \subset \bfZ, \: |\Lambda | < \infty} \:\: \mathfrak A^{S} (\Lambda )  .  $$
If $\Lambda$ is a finite set, and let $tr_{\Lambda}$ be the partial trace $tr_{\Lambda}$ satisfying
\begin{equation}
\label{eqn:a201}
 tr_{\Lambda} (QR) = Q  tr(R)  \in \bfC \cdot 1 \quad
 (\: Q \in \mathfrak A^{S}(\Lambda^{c}), \:\:  R \in \mathfrak A^{S} (\Lambda) \: )
\end{equation}
$ tr_{\Lambda} $ is a projection from  $\mathfrak A^{S}$ to $\mathfrak A^{S} (\Lambda^{c} )$  constructed in the following manner:
\begin{equation}
\label{eqn:a202}
tr_{\Lambda} (Q) =   c \sum_{k,l} \:\: e_{kl} Q e_{lk}
\end{equation}
where $\{ e_{kl} \in \mathfrak A^{S} (\Lambda ) \}$ is a matrix unit system of $\mathfrak A^{S} (\Lambda )$, 
i.e. the linear hull of $\{ e_{kl} \}$ is $\mathfrak A^{S} (\Lambda )$ and  $ e_{kl} $ satisfy
$$e_{kl}e_{ab} = \delta_{la} e_{kb} , \quad e_{kl}^{*} = e_{lk} , \quad \sum_{k} e_{kk} = 1,$$
and $c$ is the normalization constant  determined by $tr_{\Lambda} (1) = 1$.
\par  
Due to realization (\ref{eqn:a202}), $tr_{\Lambda}$ can be extended to a normal projection from  
$\pi (\mathfrak A^{S})^{\dprime}$ to $\mathfrak A^{S} (\Lambda^{c} )^{\dprime}$ 
for any representation $\pi (\cdot )$ of $\mathfrak A^{S}$.
\par
We turn to our proof of Lemma \ref{lmm:2.2} (i). 
Let $\{ \pi_{\varphi^{S}}(\mathfrak A^{S}) , \Omega_{\varphi^{S}} , \mathfrak H_{\varphi^{S}} \}$ be the GNS triple associated with 
$\mathfrak A^{S} (\Lambda )$. 
\begin{equation}
\label{eqn:a212}
\varphi (Q^{*}R^{*}RQ)=\varphi (R^{*}Q^{*}QR)\leq ||Q||^{2} \varphi (R^{*}R)
\end{equation}
for $Q \in \mathfrak A^{S}_{L},\:\: R \in \mathfrak A^{S}_{R}$,  and 
the representation  of $\pi_{\varphi^{S}}(\mathfrak A_{R}^{S})$ on 
$ \mathfrak H_{\varphi^{S}}$ and that on the closure of $ \pi_{\varphi^{S}}(\mathfrak A_{R}^{S})\Omega_{\varphi^{S}}$ are quasi-equivalent
and we can identify the center of $\pi_{\varphi^{S}}(\mathfrak A_{R}^{S})^{\dprime}$  on both spaces. 
\par 
We introduce the algebra at infinity via the following equation: 
$$\mathfrak M_{\infty} =  \cap_{ n= 1,2,3\cdots}  \:\: \pi_{\varphi^{S}}(\mathfrak A^{S} ([-n,n]^{c})^{\dprime},$$
and we claim that the center of $\pi_{\varphi^{S}}(\mathfrak A^{S})^{\dprime}$ coincides with $\mathfrak M_{\infty}$. By definition, 
inclusion $\mathfrak M_{\infty} \subset \pi_{\varphi^{S}}(\mathfrak A^{S})^{\dprime} \cap \pi_{\varphi^{S}}(\mathfrak A^{S})^{\prime}$
follows. 
If $C$ is an element of the center of $\pi_{\varphi^{S}}(\mathfrak A^{S})^{\dprime}$, there exists  
$C_{\alpha} \in \mathfrak A^{S}_{loc}$ such that  $C = w-\lim_{\alpha \to\infty} C_{\alpha}$. Then, for each $n>0$,
$$w-\lim_{\alpha \to\infty} tr_{[-n,n]} (C_{\alpha}) =  tr_{[-n,n]} (C) = C$$
which tells us 
$\mathfrak M_{\infty} \supset \pi_{\varphi^{S}}(\mathfrak A^{S})^{\dprime} \cap \pi_{\varphi^{S}}(\mathfrak A^{S})^{\prime}$.
In the same line of reasoning, the algebra at infinity of $\pi_{\varphi^{S}}(\mathfrak A_{R}^{S})^{\dprime}$  is the
center of $\pi_{\varphi^{S}}(\mathfrak A_{R}^{S})^{\dprime}$ . As any element in $\pi_{\varphi^{S}}(\mathfrak A_{R}^{S})^{\dprime}$
satisfies the condition of elements of the algebra at infinity of $\pi_{\varphi^{S}}(\mathfrak A^{S})^{\dprime}$ and 
the algebra at infinity of $\pi_{\varphi^{S}}(\mathfrak A_{R}^{S})^{\dprime}$ is contained in that of 
$\pi_{\varphi^{S}}(\mathfrak A^{S})^{\dprime}$. If $\varphi$ is a factor state of $\mathfrak A^{S} (\Lambda )$
$\pi_{\varphi^{S}}(\mathfrak A_{R}^{S})^{\dprime}$ is factor.
\par
Next, we consider the fermionic case (ii) of Lemma \ref{lmm:2.2}.  We assume that $\psi$ is a factor state of $\mathfrak A$.
As $\psi$ is $\Theta$ invariant, $\Theta$ can be extended to an automorphism of  $\pi_{\psi}(\mathfrak A )^{\dprime}$ and
any central element $C$ of $\pi_{\psi}(\mathfrak A_{R} )^{\dprime}$ is  a sum of even and odd elements $C_{\pm}$, 
$\Theta (C_{\pm})=\pm C_{\pm} , C = C_{+} + C_{-}$.
\par
If we identify the one-particle space $K$ with  $l^{2}(\bfZ )\oplus l^{2}(\bfZ )$,  $K_{L}$ with 
$l^{2}(\bfZ_{< 0})\oplus l^{2}(\bfZ_{< 0})$ and   $K_{R}$ with $l^{2}(\bfN \cup\{ 0\}) \oplus l^{2}(\bfN \cup\{ 0\})$ 
where $\bfZ_{< 0} = \{ k \in \bfZ \:|\:  k < 0 \}$.
$\mathfrak A^{S}_{L/R}$ is isomorphic to $\mathfrak A_{L/R}$. However, as  odd elements of $\mathfrak A^{S}_{L}$ and 
$\mathfrak A^{S}_{R}$  are anticommuting,  we can only apply our argument for (i) of Lemma \ref{lmm:2.2} for the even part $C_{+}$.
( c.f. Notes and Remark of Section2.6.1 of \cite{BratteliRobinsonI} )
If  $C \in \pi_{\psi}(\mathfrak A_{R} )^{\dprime}\cap \pi_{\psi}(\mathfrak A_{R} )^{\prime}$ and $\Theta (C) =C$, 
, triviality $C = d 1$  holds  for some$d \in \bfC$.
\par
If  $C =  - \Theta (C)\ne 0$ is self-adjoint, belonging to the center of  $\pi_{\psi}(\mathfrak A_{R} )^{\dprime}$,then, $\Theta (C^{2} ) = d 1$
with $d > 0$,and $X = d^{-1/2} C$ is a self-adjoint unitary. $X $ anticommutes with odd elements of $\pi_{\psi}(\mathfrak A_{L} )^{\dprime}$.
Thus, $X \in \pi_{\psi}(\mathfrak A_{R} )^{\dprime}$ is a self-adjoint unitary satisfying
\begin{equation}
X\pi_{\psi}(Q) X^{-1} =  \pi_{\psi}(\Theta_{-}(Q))   \quad  Q \in \mathfrak A.
\label{eqn:b203}
\end{equation}
In another word, $X$ implements $\Theta_{-}$.
We arrive at our claim of (ii) of Lemma \ref{lmm:2.2}. For any factor state $\psi$,  $\psi$ and $\psi \circ\Theta_{-}$ are quasi-equivalent or
disjoint. If  $\psi$ and $\psi \circ\Theta_{-}$ are disjoint, $X$ satisfying (\ref{eqn:b203}) cannot exist and 
$\pi_{\psi}(\mathfrak A_{R} )^{\dprime}$ is a factor. This implies (ii-a) of Lemma \ref{lmm:2.2}.
\par
If $\psi$ is pure,  $\psi$ restricted to $\mathfrak A^{(+)}$ and $\psi \circ\Theta_{-}$ restricted to $\mathfrak A^{(+)}$
are either  unitarily equivalent or disjoint. If they are unitarily equivalent , $\pi_{\psi}(\mathfrak A_{R} )^{\dprime}$ is a factor
because if $X$ satisfying (\ref{eqn:b203}) exists, $\psi$ cannot be pure due to Lemma \ref{lmm:2.1}.
This shows (ii-b) of Lemma \ref{lmm:2.2}.
\par
If $\psi$ is pure, and if $\psi$ and $\psi \circ\Theta_{-}$ are unitarily equivalent, and
 $\psi$ restricted to $\mathfrak A^{(+)}$ and $\psi \circ\Theta_{-}$ restricted to $\mathfrak A^{(+)}$
are disjoint, there exists a self-adjoint unitary X satisfying satisfying (\ref{eqn:b203}) again due to Lemma \ref{lmm:2.1}.
Let $\psi_{R}$ be the normal extension of $\psi$ to $\pi_{\psi}(\mathfrak A_{R} )^{\dprime}$ and
Set
$$ \psi_{1}(Q) = \psi_{R}( (1+X)Q) ,\quad   \psi_{2}(Q) = \psi_{R}( (1-X)Q) \quad Q \in \pi_{\psi}(\mathfrak A_{R} )^{\dprime} .$$
As $\psi_{R} = 1/2( \psi_{1} +\psi_{2})$, we obtain the factor decomposition of $\psi_{R}$ which implies 
(ii-c) of Lemma \ref{lmm:2.2}.
%%%%%%%%%%%%%%%%%%%
%%%%%%%%%%%%%%%%%%%
\\
\\
In Lemma \ref{lmm:2.3} we consider equivalence of projections in a von Neumann algebra.
All the knowledge we need here is explained in  Chapter 6 of \cite{Kadison}.
Let $\mathfrak M$ be a von Neumann algebra and $p,q$ be projections in $\mathfrak M$.
$p,q$ are equivalent if and only if there exists a partial isometry $u \in \mathfrak M$ satisfying $uu^{*}=p$ and $u^{*}u = q$
If $\mathfrak M$ is a factor, any two projections are comparable, i.e. either (a) $p$ is equivalent to a sub-projection of $q$ 
or (b) $q$ is equivalent to a sub-projection of $p$.  
If $\mathfrak M$ is a type III factor, any two projections are equivalent. If $\mathfrak M$ is a type II factor with a normal semi-finite trace $tr$, 
 two projections $p,q$  are equivalent if and only if $tr(p) = tr(q)$.
\begin{lmm}
\label{lmm:2.3}
Let $\mathfrak M$ be a von Neumann algebra with a finite center 
acting on a Hilbert space $\mathfrak H$ and $\Gamma$ be a self-adjoint unitary on $\mathfrak H$
giving rise to a $Z_{2}$ action $\Theta$ on $\mathfrak M$,
and set $\mathfrak N = \{ Q \in \mathfrak M \:| \: \Gamma Q = Q\Gamma \: \}$. If 
$\mathfrak N$ is of type I with a finite center, $\mathfrak M$ is a type I von Neumann algebra.
\end{lmm}
{\it Proof of Lemma \ref{lmm:2.3}}
\\
Due to our assumption of Lemma \ref{lmm:2.3}
, $\mathfrak M$ cannot be a finite factor , we proceed to our proof assuming that $\mathfrak M$ is either of type $II_{\infty}$ or of type $III$.
\par 
Suppose that $\mathfrak M$ is of type $II_{\infty}$ and let $tr(\cdot )$ be a $\Theta$ invariant semi-finite normal trace of $\mathfrak M$.
Let $P$ be a minimal projection in   a factor component of $\mathfrak N$. For any small  positive number $c$ satisfying 
$0 < c < 1/2 d = 1/2 tr(P)$ 
there exists a projection $Q \in \mathfrak M$ such that $ Q \leq P , c=tr(Q)$. 
Set $Q_{\pm} =1/2 (Q\pm\Theta (Q))$. As  $tr(\cdot )$ is $\Theta$ invariant, $c= tr(Q_{+})$
\par 
Suppose $d = \infty$. As $Q_{+}  \leq P, Q_{+} \in \mathfrak N$ and  $P$ is minimal, $Q_{+} = a P$ for some $a > 0$ which leads to
a contradiction , $ c = tr(Q) = a tr(P) = \infty$. Suppose $d$ is finite. 
As $Q$ is a projection, 
$$Q^{2} = Q_{+}^{2} + (Q_{+}Q_{-}+ Q_{-}Q_{+})+ Q_{-}^{2} = Q \leq P .$$ 
Applying $\Theta$,
we obtain 
\begin{equation}
\label{eqn:b201}
0\leq  Q_{+} = Q_{+}^{2} + Q_{-}^{2} \leq P , \quad Q_{-} = Q_{+}Q_{-}+ Q_{-}Q_{+}.
\end{equation}
As  $Q_{+} $, and $Q_{-}^{2}$  are positive elements in $\mathfrak N$ 
and  $P$ is a minimal projection, we obtain $Q_{+} =cP, Q_{-}^{2} = (c- c^{2}) P$.
Then,
$$PQ_{-}P =  PQ_{+}Q_{-}P + PQ_{-}Q_{+}P = 2c PQ_{-}P $$
and as a consequence, $PQ_{-}P = 0$ and 
Due to the second identity of (\ref{eqn:b201}),
$$PQ_{-} = PQ_{+}Q_{-}+ PQ_{-}Q_{+} = cPQ_{-}+ cPQ_{-}P = cPQ_{-} .$$
This shows that $PQ_{-} = 0$ and  that $Q_{-} = Q_{+}Q_{-}+ Q_{-}Q_{+}= c PQ_{-} +c Q_{-}P = 0$.
It turns out   $Q = Q_{+} + Q_{-} = Q_{+} \in \mathfrak N$, and $P= Q$ 
which leads to a contradiction.
\par
Next we assume that $\mathfrak M$ is a type $III$ von Neumann algebra. We can take a $\Theta$ invariant central projection
$P \ne 0 $ of $\mathfrak M$ and a projection Q of a rank, at most 2,  in  in one factor component of $\mathfrak N$ . 
( If $P_{1}$ is a minimal central projection of $\mathfrak M$ such that $\Theta (P_{1}) P_{1}= 0$ we set $P = P_{1} + \Theta (P_{1})$ 
and we take  a minimal projection $Q_{1}$ of $\mathfrak N$ equivalent to $P_{1}$ and a minimal projection $Q_{2}$ of $\mathfrak N$ 
equivalent to $P_{2}$ orthogonal to $Q_{1}$.
Then , $P$ is equivalent to $Q= Q_{1}+Q_{2}$.) We now consider $\mathfrak M P$ and by abuse of notation $P =1$, $\mathfrak M=\mathfrak M P$
\par
Let $U$ be a partial isometry satisfying $Q = UU^{*} , 1 = U^{*}U$. We can write $U = U_{+} + U_{-}X$ where $X =X^{*} = - \Theta (X), X^{2}=1$, $U_{\pm} \in \mathfrak N$. As $Q, 1 \in \mathfrak N$ ,
$$Q = U_{+}U_{+}^{*} +  U_{-}U_{-}^{*}, \quad 1  = U_{+}^{*}U_{+} +  XU_{-}^{*} U_{-}X  $$
$$U_{+}^{*} U_{-}X + XU_{-}^{*} U_{+}=0, \quad U_{-}XU_{+}^{*}  + U_{+}XU_{-}^{*} =0,$$
As  $ U_{\pm}U_{\pm}^{*}\leq Q$, $U_{\pm}= QU_{\pm}$ which means the rank( the dimension of the range) of $U_{\pm}$ is at most two.
As a consequence, the rank of  $U_{+}^{*}U_{+}$ and that of $XU_{-}^{*} U_{-}X$ are at most two and 
$U_{+}^{*}U_{+} +  XU_{-}^{*} U_{-}X$ is of finite rank , which is a contradiction.
\bigskip
\newline
{\it Proof of Theorem \ref{th:1.3}}
\\
Let $\psi$ be a $\Theta$ invariant pure state of $\mathfrak A$ and we focus on the restriction of $\psi$ to $\mathfrak A_{L}$.
We fix  $h \in K_{L}$ such that  $Jh=h$ and $|| h || =1$.  Then,  $B(h)$ is a self-adjoint unitary in $\mathfrak A_{L}$.
Let $\{ \pi_{\psi}(\mathfrak A ), \mathfrak H_{\psi} , \Omega_{\psi} \}$  be  the GNS triple associated with $\psi$.    
Set $\mathfrak H_{\psi}^{(\pm)} =  \overline{ \pi (\mathfrak A^{(\pm )} )\Omega_{\psi}}$ and by $\pi_{\pm} (\mathfrak A^{(+)})$ 
we denote  the representation  $\pi_{\psi} (\mathfrak A^{(+)})$ of   $\mathfrak A^{(+)}$ restricted to $\mathfrak H_{\pm}$. 
As we mentioned before, $\pi_{\psi} (\mathfrak A^{(+)})$ are disjoint.
By definition $\pi_{\psi} (\mathfrak A^{(+)})$ of   $\mathfrak A^{(+)}$ on $\mathfrak H_{\psi}$  is a direct sum of 
$\pi_{\psi} (\mathfrak A^{(+)})$. 
As $\mathfrak H_{\psi}^{(\pm)} =  \overline{ \pi (\mathfrak A^{(+)} B(h))\Omega_{\psi}}$,  
 $\pi_{\psi} (\mathfrak A^{(+)})$ of   $\mathfrak A^{(+)}$ on $\mathfrak H_{\psi}$ is unitarily equivalent to 
 the GNS representation associated with $1/2 \{ \psi + \psi\circ Ad(B(h) \}$.
\par
Suppose that $\psi$ is quasi-equivalent to $\psi_{L}\otimes\psi_{R}$ and that both $\psi_{R}$ and $\psi_{L}$ are $\Theta$ invariant.
As  to the restriction to $\mathfrak A^{(+)}$ the we have three possible case.
\\
(Case-1) $\psi$ restricted to $\mathfrak A^{(+)}$ is quasi-equivalent to $\psi_{L}\otimes\psi_{R}$ restricted to $\mathfrak A^{(+)}$.
\\
(Case-2)  $\psi\circ Ad(B(h)$ restricted to $\mathfrak A^{(+)}$ is quasi-equivalent to $\psi_{L}) \otimes\psi_{R}$ restricted to $\mathfrak A^{(+)}$.
\\
(Case-3) $1/2 \{ \psi + \psi\circ Ad(B(h) \}$  restricted to $\mathfrak A^{(+)}$ is quasi-equivalent to $\psi_{L}\otimes\psi_{R}$ restricted to $\mathfrak A^{(+)}$.
\\
Note that in all cases,  the state  $\psi$ , $\psi\circ Ad(B(h)$ or $1/2 \{ \psi + \psi\circ Ad(B(h) \}$
restricted to $\mathfrak A^{(+)}$ extended to a $\Theta$ invariant state of $\mathfrak A^{S}$ is of type I.
\par
Let $\varphi$ be the $\Theta$ invariant extension of  $\mathfrak A^{(+)}$ to $\mathfrak A^{S}$ and 
let $\varphi_{L}$ be the $\Theta$ invariant extension  of  $\psi_{L}$ restricted to $\mathfrak A_{L}^{(+)}$ to 
$\mathfrak A_{L}^{S}$ and   $\varphi_{R}$ be the $\Theta$ invariant extension  of  $\psi_{R}$ restricted to $\mathfrak A_{R}^{(+)}$ 
to $\mathfrak A_{L}^{S}$. 
\par
If $Q \in \mathfrak A_{L}$ $ R\in \mathfrak A_{R}$ satisfying $QR \in \mathfrak A^{(+)}$ we have $QR = Q_{+}R_{+} +Q_{-}R_{-}$
and
\begin{eqnarray}
\label{eqn:b3}
&& \psi_{L}\otimes\psi_{R}(QR) =  \psi_{L}(Q_{+}\psi_{R} (R_{+}) =  \varphi_{L}(Q_{+}) \varphi_{R}(R_{+}) 
\nonumber
\\
&&=\varphi_{L}\otimes\varphi_{R}((Q_{+}+ Q_{-}T ) (R_{+}+R_{-}T) ).
 \end{eqnarray}
(\ref{eqn:b3}) shows that if we extend  $\psi_{L}\otimes\psi_{R}$ restricted to $\mathfrak A^{(+)}$ to a $\Theta$ invariant state of 
$\mathfrak A^{S}$, we obtain  a product state $ \varphi_{L}\otimes\varphi_{R}$. 
If  (Case-1) is valid, $\varphi$ is quasi-equivalent to $\varphi_{L}\otimes\varphi_{R}$ due to Lemma \ref{lmm:2.0}.  
and  $\varphi$ restricted to $\mathfrak A_{L}^{S}$
which shows $\psi$ restricted to $\mathfrak A_{L}$ is of type I.
The same argument is valid for (Case-2) and (Case-3).
%%%%%%%%%%%%%%%%%%%%%%%%%
%%%%%%%%%%%%%%%%%%%%%%%%
\section{Fermionic String Order}\label{String}
\setcounter{theorem}{0}
\setcounter{equation}{0}
First we begin with Fock states of our CAR algebra  $\mathfrak A$.
A Bogoliubov automorphism $\beta_{u}$ is  an automorphism of   $\mathfrak A$
defined by the following equation:
$$ \beta_{u}(B(f)) = B(uf) \quad f \in K$$ 
where $u$ is a unitary on $K$ satisfying $JuJ = u$.
On the one-particle $K = K_{R}\oplus K_{L}$, we introduce self-adjoint unitaries $\theta$ and $\theta_{-}$ satisfying the following equations:
\begin{eqnarray}
\label{eqn:c1} 
&&  \theta (f) = - f \quad f \in K ,
\nonumber
\\
&& \theta_{-}(f_{L}\oplus f_{R}) = (- f_{L}\oplus f_{R}),\quad  f_{L} \in K_{L} ,\:\: f_{R} \in K_{R}.
\end{eqnarray}
Obviously, $\Theta (Q) = \beta_{\theta} (Q)$, $\Theta_{-}(Q) = \beta_{\theta_{-}} (Q)$ for $Q \in \mathfrak A$.
%% definition of basis projection %%%%%%
\begin{df}
\label{df:3.1}
(i) A projection $E$ acting on $K$ is a basis projection if  and only if $JEJ = 1-E$.
\\
(ii) Let $E$ be a basis projection. Let $\psi_{E}$ be a pure state  of $\mathfrak A (K)$ uniquely determined by the following equation:
\begin{equation}
\label{eqn:c2} 
\psi_{E} (B^{*}(f) B(g) ) = (f,Eg)_{K} \quad f,g \in K .
\end{equation}
We call $\psi_{E}$ a Fock state associated with $E$. 
\end{df}
 In the above Fock state,   $\psi_{E} (B^{*}((1-E)g) B((1-E)g) ) = 0$ implies 
  $B((1-E)g)$ plays a role of an annihilation operator and  $B((1-E)g)^{*} = B(J(1-E)g) = B(EJg)$ plays a role of  a creation operator. The GNS Hilbert space is an antisymmetric Fock space with a one-particle space $EK$.   
\par
 It is known that two Fock states $\psi_{E_{1}}$ and $\psi_{E_{2}}$ are unitarily equivalent if and only if $E_{1} - E_{2}$ is a Hilbert-Schmidt operator on $K$. (See \cite{Araki1971}.) In particular, for a unitary $u$ satisfying $JuJ =u$ and for a basis projection, 
if $uEu-E$ is a Hilbert-Schmidt operator, there exists a unitary $\Gamma (u)$ on the GNS Hilbert space implementing $u$: 
$$  \Gamma (u) \pi_{\psi_{E}}(Q) \Gamma (u)^{*} =   \pi_{\psi_{E}}(\beta_{u}(Q)) ,\quad Q \in \mathfrak A .$$
If  $E_{1} - E_{2}$ is a compact operator, the dimension $dim \: E_{1}\wedge (1- E_{2})$ of the range of $E_{1}\wedge (1-E_{2})$ is finite.
\begin{theorem}
\label{th:3.2}
(i) Let $E$ be a basis projection. The  split property is valid for  $\psi_{E}$ if and only if
$ \theta_{-} E \theta_{-} - E$ is a Hilbert-Schmidt operator on $K$.
\\
(ii)  If  $\psi_{E}$ has the split property, the restriction of $\psi_{E}$ to $\mathfrak A_{L}$ is factor if and only if
$dim \: E \wedge (1-\theta E\theta)$ is even.
\end{theorem}
As all elements for proof of Theorem \ref{th:3.2} is contained in \cite{Araki1971} and \cite{ArakiMatsui1985}, we sketch our proof briefly.
Suppose that $\psi_{E}$ is quasi-equivalent to $\psi^{1} \otimes_{\bfZ_{2}} \psi^{2}$. As quasi-equivalence classes of representations matters 
we may replace $\psi^{1}$  with another state of $\mathfrak A_{L}$ quasi-equivalent to  $\psi^{1}$ in $\psi^{1} \otimes \psi^{2}$ and
in view of the fermion version of (\ref{eqn:a212}), we may assume that $\psi^{1}$ is restriction of  $\psi_{E}$  to $\mathfrak A_{L}$
and that $\psi^{2}$ is restriction of  $\psi_{E}$  to $\mathfrak A_{R}$.
As a consequence, $\psi^{1} \otimes_{\bfZ_{2}} \psi^{2}$ is a quasifree state defined in Definition 3.1 of  \cite{Araki1971}
where the covariance operator $S$ in (3.3) of  \cite{Araki1971} is written as follows:
\begin{equation}
\label{eqn:c3} 
S =  qEq+ (1-q) E (1-q) , \quad q = \frac{1}{2} (1-\theta_{-}) .
\end{equation}
Due to Theorem 1 and Lemma 5.1 of  \cite{Araki1971},
$\psi^{1} \otimes \psi^{2}$ is quasi-equivalent to $\psi_{E}$ if and only if $P_{E} - P_{S}$ is in the Hilbert Schmidt class where
$P_{S}$ is defined in (4.4) of  \cite{Araki1971}. This condition is equivalent to the conditions that 
$ \theta_{-} E \theta_{-} - E$  is in the Hilbert Schmidt class and that $ \theta_{-} E \theta_{-} (1- E) \theta_{-} E \theta_{-}$ is in
the trace class, however,  the former condition implies the latter. 
(ii) of  Theorem \ref{th:3.2} is same  as Theorem 4 of \cite{ArakiEvans1983}.
\bigskip
\\
 $(-1)^{dim \: E_{1}\wedge (1- E_{2}) }$ is nothing but the classical $\bfZ_{2}$ index of real Fredholm operators and 
 based on this observation, we introduce the following  $\bfZ_{2}$ index.
\begin{df}
\label{df:3.3}
Let $\psi$ be a $\Theta$ invariant pure state of $\mathfrak A$ with split property. 
\\
(i) $ind_{\bfZ_{2}} \psi =1$  if the restriction of  $\psi$ to $\mathfrak A_{L}$
is a factor state. 
\\
(ii) $ind_{\bfZ_{2}} \psi = -1$  if the restriction of  $\psi$ to $\mathfrak A_{L}$ is not a factor state. 
\end{df}
\begin{rmk}
\label{rmk:after3.3}
 $\psi$ and  $\psi\circ\Theta_{-}$ are unitarily equivalent if split property holds for any $\Theta$ invariant pure state $\psi$  of $\mathfrak A$. 
We are wondering if the converse holds.
Namely, unitary equivalence of $\psi$ and  $\psi\circ\Theta_{-}$ implies split property or not. We are unable to prove the claim , 
unable to provide counter any example.
\end{rmk}
For automorphisms $\alpha_{L}$ of $\mathfrak A_{L}$ and $\alpha_{R}$ of $\mathfrak A_{R}$
commuting  with $\Theta$, we can introduce a $\bfZ_{2}$ graded product automorphism 
$\alpha_{L}\otimes_{\bfZ_{2}} \alpha_{R}$ of  $\mathfrak A$
satisfying 
\begin{equation}
\label{eqn:c4} 
\alpha_{L}\otimes_{\bfZ_{2}} \alpha_{R}(Q_{L}) = \alpha_{L}(Q_{L}) ,\:\:\: \alpha_{L}\otimes_{\bfZ_{2}} \alpha_{R}(Q_{R}) = \alpha_{L}(Q_{R})
\end{equation}
for $Q_{L} \in \mathfrak A_{L} ,\:\:\:  Q_{R} \in \mathfrak A_{R}$ and
\begin{equation}
\label{eqn:c5}
(\varphi_{L}\otimes_{\bfZ_{2}} \varphi_{R}) \circ (\alpha_{L}\otimes_{\bfZ_{2}} \alpha_{R}) 
= (\varphi_{L}\circ \alpha_{L}) \otimes_{\bfZ_{2}} (\varphi_{R}\circ\alpha_{R})
\end{equation}
for any $\bfZ_{2}$ graded product state $\varphi_{L}\otimes_{\bfZ_{2}} \varphi_{R}$.
\begin{df}
\label{df:3.4}
An automorphism $\alpha$ of  $\mathfrak A$ is an almost product automorphism if there exists a 
$\bfZ_{2}$ graded product automorphism $\alpha_{L}\otimes_{\bfZ_{2}} \alpha_{R}$  and a unitary $U$ of  $\mathfrak A$
such that
\begin{equation}
\label{eqn:c6} 
\alpha   = ( \alpha_{L}\otimes \alpha_{R} )\circ Ad(U).
\end{equation}
\end{df}
As an inner perturbation $Ad(U)$ preserves quasi-equivalence classes of representations, the following proposition
is obvious.
\begin{pro}
\label{pro:3.5}
Let $\psi$ be a $\Theta$ invariant pure state of  $\mathfrak A$ and let $\alpha$ be an almost product automorphism
of $\mathfrak A$. Then,
\begin{equation}
\label{eqn:c7} 
ind_{\bfZ_{2}} \psi = ind_{\bfZ_{2}} (\psi\circ\alpha )
\end{equation}
\end{pro}
Next we look into gapped ground states of one-dimensional Fermion systems.   
The one-particle space we deal with is $K = l^{2}(\bfZ ) \oplus l^{2}(\bfZ ) $. 
We assume that systems are translationally invariant  without loss of generality.  
Our argument below is valid for periodic systems  or for one-particle subspace with an internal degree of freedom $V$, 
$$K =  l^{2}(\bfZ )\otimes V \oplus   l^{2}(\bfZ )\otimes V$$ 
provided the dimension of $V$ is finite.
The correspondence of creation and annihilation operators $c^{*}_{k}$ and $c_{l}$
and $B(f\oplus g)$ is as follows. 
\par
Set $\bfZ_{<} = \{ j \in \bfZ \: | \: j < 0 \}$,  $K= l^{2}(\bfZ )\oplus l^{2}(\bfZ )$,
 $ K_{L}= l^{2}(\bfZ_{<} ) \oplus  l^{2}(\bfZ_{<} )$, $K_{R}= l^{2}(\bfN \cup \{0 \} ) \oplus  l^{2}(\bfN \cup \{0 \} ) $
\begin{eqnarray}
\label{eqn:c702} 
&&B(f \oplus g ) = c^{*}(f) + c(g) =  \sum_{j} \left( c^{*} f_{j} + c g_{j}\right)
\nonumber
\\
&&J(f \oplus g ) = ( \overline{g} \oplus \overline{f} )
\end{eqnarray}
for $f = \{f_{j}\} \in K_{0} , g = \{g_{j}\} \in K_{0}$ ,  $\overline{f}$ is the complex conjugation $\overline{f} =  \{\overline{f}_{j}\}$.
\begin{pro}
\label{pro:3.6}
Let $\psi$ be a translationally invariant pure ground state of  a finite range translationally invariant  Hamiltonian of (\ref{eqn:a13-2}).
If $\psi$ is a gapped ground state, the split property holds, i.e. $\psi$ restricted to $\mathfrak A_{L}$ is type I. 
\end{pro}
We may prove Proposition \ref{pro:3.6} in several ways.
As repeating proof of \cite{Hastings2007}, \cite{Brandao2015 }  and \cite{Matsui2010} for fermion systems is somewhat lengthy,
here we employ the Jordan-Wigner transformation a la mani$\grave{e}$re de \cite{ArakiMatsui1985}
and reduce the problem to the Pauli spin chain. We work in the algebra $\hat{\mathfrak A}$ and $T$ to define
our Jordan-Wigner transformation. For any $l \in \bfZ$, set
\begin{eqnarray}
\label{eqn:c8}
&&\sigma_{x}^{(l)} = S_{l} (c_{l}^{*} + c_{l} )
\nonumber
\\
&&\sigma_{y}^{(l)} =  i S_{l} (c_{l} - c_{l}^{*} )
\nonumber
\\
&&\sigma_{z}^{(l)} = ( 2c^{*}_{l}c_{l} -1 )
\end{eqnarray}
where
$$
S_{j} = \begin{cases} \quad   T \prod_{k=o}^{j-1}\sigma_{z}^{(k)}  \quad ( 1 \leq j )
\\    
\quad\quad\quad\:\: T   \quad\quad\quad ( j = 0 )
\\    
\quad T  \prod_{k=j}^{-1}\sigma_{z}^{(k)}  \quad ( j <  0 ).
\end{cases}
$$
We obtain relations for Pauli spin matrices:
$$ \sigma_{x}^{(l)} \sigma_{y}^{(l)} = i  \sigma_{z}^{(l)} , \quad \{\: \sigma_{x}^{(l)}\: , \: \sigma_{y}^{(l)} \} = 0,\quad   
(\sigma_{\alpha}^{(l)})^{2} = 1  \quad [ \sigma_{\alpha}^{(l)} , \sigma_{\beta}^{(k)} ] = 0$$
for $l , k \in \bfZ,  l \ne k$ and   $\alpha, \beta = x,y,z$.
\par
Let $\mathfrak A^{P}$ be a $C^{*}$-subalgebra of $\hat{\mathfrak A}$ generated by $\sigma_{\alpha}^{(l)} \:\: l \in\bfZ, \:\: \alpha = x,y,z$
where $P$ stands for {\it Pauli spin algebra}. 
Set
$$\mathfrak A^{P \; (\pm )} = \{ Q \in \mathfrak A^{P} \: | \: \Theta (Q) = \pm Q \}, 
\quad Q_{\pm} = \frac{Q \pm \Theta (Q)}{2} \in \mathfrak A^{P \; (\pm )}$$ 
and let $\mathfrak A_{[n,m]}^{P}$ be the  $C^{*}$-subalgebra of $\mathfrak A^{P}$ generated by   
$\sigma_{\alpha}^{(l)} \:\:n \leq l \leq m, \:\: \alpha = x,y,z$. 
Set 
$$\mathfrak A_{[n,m]}^{P (\pm)}  =   \mathfrak A_{[n,m]}^{P} \cap \mathfrak A^{P \; (\pm )} , $$
$$\mathfrak A_{loc}^{P}  \:\:\: = \:\:\: \cup_{-\infty < n \leq m <\infty} \mathfrak A_{[n,m]} .$$
, i.e. $\mathfrak A_{loc}^{P}$ is the set of all elements strictly localized in finite sets.  
Finally  we introduce $\mathfrak A_{L}^{P} =  \mathfrak A_{(-\infty ,-1]}^{P} $ , $\mathfrak A_{R}^{P} =  \mathfrak A_{[0, \infty )}^{P} $.
By split property of a state $\varphi$ of $\mathfrak A^{P}$ we mean split property with respect to
 $\mathfrak A_{L}^{P}$ and $\mathfrak A_{R}^{P}$.
\par
Due to our definition, it easy to see that $\mathfrak A^{(+)}$ and  $\mathfrak A^{P \; (+)}$
coincide: $\mathfrak A^{(+)} = \mathfrak A^{P \; (+)}$. 
The lattice translation $\tau_{j}$ on $\mathfrak A^{(+)}$ can be extended to $\mathfrak A^{P}$.
\par
For any $\Theta$ invariant pure state $\psi$ of $\mathfrak A$, we extend $\psi$ to  a state $\hat{\psi}$ of $\hat{\mathfrak A}$
via the following equation:
$$ \hat{\psi} (Q_{1} + Q_{2} T ) = \psi (Q_{2} ) \quad Q_{1} , Q_{2} \in \mathfrak A .$$
$\hat{\psi}$ is invariant under $\Theta$.
\par
We recall results proved in \cite{ArakiMatsui1985}. See Section 5 of \cite{ArakiMatsui1985}.  In our proof Section 5 of \cite{ArakiMatsui1985},
quasifree property of states is not used but a crucial point is (non-)existence of $U(\Theta_{-})$ 
where  $U(\Theta_{-})$ is the self-adjoint unitary implementing $\Theta_{-}$ on $\mathfrak A$ on the GNS space of $\psi$.  
Note that if  $\psi$ be a $\Theta$ invariant pure state of $\mathfrak A$ and if $\psi$ and  $\psi\circ\Theta_{-}$ are unitarily equivalent,
there exists a self-adjoint unitary $U(\Theta_{-})$  on the GNS space $\mathfrak H_{\psi}$ satisfying
$$U(\Theta_{-}) \pi_{\psi}((Q)) U(\Theta_{-})^{*} = \pi_{\psi}(\Theta_{-}(Q)) ,\quad Q \in \mathfrak A.$$
If we assume further that $\psi$ and  $\psi\circ\Theta_{-}$ restricted to $\mathfrak A^{(+)}$ are unitarily equivalent,
$U(\Theta_{-}) \in \pi_{\psi}(\mathfrak A^{(+)})^{\dprime}$. If we assume, instead,  that
$\psi$ and  $\psi\circ\Theta_{-}$ restricted to $\mathfrak A^{(+)}$ are disjoint, 
 $U(\Theta_{-})$ is in the closure of $\pi_{\psi}(\mathfrak A^{(+)})^{\dprime}$ by weak operator topology.
\begin{pro}
\label{pro:3.7}
Let $\psi$ be a $\Theta$ invariant pure state of $\mathfrak A$ and $\psi^{P}$ be the state
of $\mathfrak A^{P}$ uniquely determined by the following equation :
\begin{equation}
\label{eqn:c21}
\psi^{P} (Q)  = \psi (Q_{+}), \quad Q \in \mathfrak A^{P}.
\end{equation}
(i) If $\psi$ and  $\psi\circ\Theta_{-}$ are disjoint,  $\psi^{P}$ is a pure state of $\mathfrak A^{P}$. The split property does not hold.
\\
(ii) Suppose $\psi$ and  $\psi\circ\Theta_{-}$ are unitarily equivalent. 
\\
(ii-a) Assume further that $\psi$ and  $\psi\circ\Theta_{-}$ restricted to $\mathfrak A^{(+)}$ are unitarily equivalent, 
$\psi^{P}$ is a pure state of $\mathfrak A^{P}$. 
\\
(ii-b)  Assume further that  $\psi$ and  $\psi\circ\Theta_{-}$ restricted to $\mathfrak A^{(+)}$ are disjoint, 
$\psi^{S}$ is not a pure state of $\mathfrak A^{P}$.
There exists a pure state $\varphi$ of $\mathfrak A^{P}$ satisfying
\begin{equation}
\label{eqn:c22}
\psi^{P} = 1/2 \left( \varphi + \varphi\circ\Theta  \right), \quad \varphi  \ne\varphi\circ\Theta .
\end{equation}
$\varphi$ and  $\varphi\circ\Theta$ are disjoint and $\varphi$ is described as
\begin{equation}
\label{eqn:c22b}
 \varphi (Q) = \psi ( (1+ U(\Theta_{-})T) Q) ,\quad Q \in \mathfrak A^{P} 
 \end{equation}
 where by abuse of notations, we denote the normal extension of $\varphi$ to $\pi_{\psi^{P}}(\mathfrak A^{P})^{\prime}$
 by the same symbol $\varphi$. $U(\Theta_{-})$ is the self-adjoint unitary implementing $\Theta_{-}$ on 
 $\pi_{\psi}(\mathfrak A)^{\dprime}$ and $U(\Theta_{-})$ is in the weak limit of elements in $\pi_{\psi}(\mathfrak A^{(-)})$. 
 The existence of $U(\Theta_{-})$ follows from disjointness of  $\psi$ and  $\psi\circ\Theta_{-}$ restricted to $\mathfrak A^{(+)}$.
 \end{pro}
 As any $\bfZ_{2}$ product state is $\Theta_{-}$ invariant , for a $\Theta$ invariant pure state $\psi$ of $\mathfrak A$ with split property
 $\psi$ and  $\psi\circ\Theta_{-}$ are unitarily equivalent.  Proposition \ref{pro:3.7} and (ii) of Lemma \ref{lmm:2.2}  imply the following. 
\begin{cor}
\label{cor:3.7b}
 Let $\psi$ be a $\Theta$ invariant pure state of $\mathfrak A$  with split property.
\\
(i) if $ind_{\bfZ_{2}} \psi =1 $, $\psi^{P}$ is a pure state.
\\
(ii) if  if $ind_{\bfZ_{2}} \psi = -1 $, $\psi^{P}$ is not a pure state and  a pure state $ \varphi$ satisfying (\ref{eqn:c22}) and (\ref{eqn:c22b})
exists.
\\
Furthermore, if $\psi$ is periodic with a period $p$, $\psi\circ\tau_{p}=\psi$,  $\psi^{P}$ is periodic with a period $p$. 
In the case of (ii), one of the following two possibility (a), (b) is valid.
\\
(a) $\varphi$ is periodic with a period $p$.
(b) $\varphi\circ\tau_{p}=\varphi\circ\Theta$.
\end{cor}
%%%%%%% '±'±'©'ç
%%%%%%%
%%%%%%%

 It is now time to talk about gapped ground states of $\mathfrak A$. The finite volume Hamiltonian   $H_{M}$ 
 of (\ref{eqn:a13}) belongs to $\mathfrak A^{P (+)} = \mathfrak A^{(+)}$ and generates an infinite volume time evolution
 of  $\mathfrak A^{P}$ via the same limit as in (\ref{eqn:a14}). 
 We denote the time evolution of $\mathfrak A^{P}$ by the same symbol $\evl^{h}$:
\begin{equation}
\label{eqn:c23}
 \lim_{M\to\infty} \:\: e^{it H_{M} } Q  e^{-it H_{M}}  =  \evl^{h} (Q) \quad Q \in \mathfrak A^{P} .
\end{equation}
Suppose that $\psi$ is a translationally invariant pure ground state of $\evl^{h}$ for $\mathfrak A$. 
In view of (\ref{eqn:a14}) for $\mathfrak A^{P}$
the $\Theta$ invariant extension $\psi^{P}$ in (\ref{eqn:c21}) to $\mathfrak A^{P}$ is 
a translationally invariant $\Theta$ invariant ground state 
for $\evl^{P,h}$. Due to Proposition {pro:3.7},  either (i) $\psi^{P}$ is pure or 
(ii) $\psi^{P}$ is an average of two pure states  $\varphi$ and $\varphi\circ\Theta$.
\begin{pro}
\label{pro:3.8}
Suppose that $\psi$ is a translationally invariant pure gapped  ground state of $\evl^{h}$ for $\mathfrak A$.
\\
(i) if $\psi^{P}$ is a pure state of $\mathfrak A^{P}$,  $\psi^{P}$ is a translationally invariant gapped ground state.
\\
(ii)  If $\psi^{P}$ is not pure, $\varphi$ in (\ref{eqn:c22})  is a gapped ground state which is periodic with period 2 .
\end{pro}
{\it Proof  of Proposition \ref{pro:3.8}}
\\
 Let $\{ \pi (\mathfrak A^{P} ), \Omega , \mathfrak H \}$ be the GNS triple associated with  $\psi^{P}$ and set
$$\mathfrak H^{\pm} =  \overline{\pi (\mathfrak A^{P (\pm )})\Omega} . $$
Let $H_{\psi^{P}}$ be a positive self-adjoint operator
on $\mathfrak H$ implementing  $\evl^{h}$ on $\pi (\mathfrak A^{P} )$.
\par
Now we show (i).
As restriction of $\psi^{P}$ to $\mathfrak A^{P (+)}$ is a gapped ground state, $H_{\psi^{P}}$ has positive spectrum on $\mathfrak H^{+}$
$$ spec \:\:  H_{\psi^{P}} |_{\mathfrak H^{+}} \subset \{ 0 \} \cup [m ,\infty )$$
Suppose that there exists a gapless  excitation on $\mathfrak H^{-}$, in another word, for any small $\delta > 0$ and some $\eta > 0$
$$ spec \:\:  H_{\psi^{P}} |_{\mathfrak H^{-}} \cap (- \eta ,\delta ) \ne \emptyset.$$
there exists $Q \in \mathfrak A^{P (-)}_{loc}$ and a smooth function $f(t)$ such that its fourier transform $\hat{f}(\xi )$ is supported in   
$(-\delta, \eta )$ and that
\begin{equation}
\label{eqn:c23b}
\pi ( \alpha_{f}^{h}(Q))\Omega = \hat{f}(-H_{\psi^{P}} )\pi (Q)\Omega \:\: \ne 0
\end{equation}
where $ \alpha_{f}^{h}(Q)$ is an operator creating small energy excitation:
$$ \alpha_{f}^{h}(Q) = \int_{\bfR} \: \alpha_{t}^{h}(Q) \: f(t) \: dt \:\: \in  \mathfrak A^{P (-)}.$$
Due to cluster property of any pure state,
\begin{eqnarray}
\label{eqn:c24}
 && \lim_{k\to\infty} \:\: \psi^{P} ( \alpha_{f}^{h}(Q)^{*} \tau_{k} ( \alpha_{f}^{h}(Q)^{*} \alpha_{f}^{h}(Q)) \alpha_{f}^{h}(Q))
\nonumber
\\
&&=  \psi^{P} ( \alpha_{f}^{h}(Q)^{*} \alpha_{f}^{h}(Q) ) \psi^{P} ( \alpha_{f}^{h}(Q)^{*} \alpha_{f}^{h}(Q) ) \ne 0.
\end{eqnarray}
As a consequence we see that , for a large $k$, 
\begin{equation}
\label{eqn:c25}
\xi = \pi ( \alpha_{f}^{P,h}(Q)) \tau_{k} (\pi ( \alpha_{f}^{P,h}(Q)) )\Omega \ne 0 .
\end{equation}
By the proof of Lemma 3.2.42 (iii) in \cite{BratteliRobinsonI},
$\xi$ has the spectral support in $(0 ,2\delta )$ in the special decomposition of  $H_{\psi^{P}}$.
$\xi$ belongs to  $\mathfrak H^{+}$ and we obtain a contradiction.
\\
(ii) If $\psi^{P}$ is not pure, representation of  $\mathfrak A^{P (+)}$ on $\mathfrak H^{\pm}$ are unitary equivalent
the spectrum of  $H_{\psi^{P}}$ on each space is same with doubled multiplicity.
If we restrict $H_{\psi^{P}}$ to irreducible representations of $\mathfrak A^{P}$, the degeneracy of the ground state energy spectrum
is removed.
%\begin{rmk}
%Even if the states are not translationally invariant, we can still define ground states , however the non-translationally invariant ground state
%of a Pauli spin system, if it exists, may not yield a ground state of the Fermion system.   
%\end{rmk}
\bigskip
\\
{\it Proof  of Proposition \ref{pro:3.6}}
\\
(i) If $\psi^{P}$ is pure, $\psi^{P}$ is a gapped ground state for the Pauli spin system and split property is valid
due to a result in \cite{Hastings2007} and \cite{Matsui2010}.  Thus $\psi_{P}(\mathfrak  A^{P}_{L})^{\dprime}$ is type I
and $\Theta$ invariance of $\psi^{P}$ restricted to $\mathfrak  A^{P}_{L}$ implies $\pi_{P}(\mathfrak  A^{P (+)}_{L})^{\dprime}$
is type I as the commutant of a single self-adjoint unitary in $B(\mathfrak H )$ is type I.
Since $\pi_{P}(\mathfrak  A^{P (+)}_{L})^{\dprime} = \pi_{\psi}(\mathfrak  A^{ (+)}_{L})^{\dprime}$ , 
$\pi_{\psi}(\mathfrak  A_{L})^{\dprime}$ is type I.
\\
(ii) If $\psi^{P}$ is not pure, $\varphi$ is a gapped ground state for the Pauli spin system and split property is valid for
$\varphi$ and $\varphi\circ\Theta$. As $\varphi$ is pure,  $\varphi$ is equivalent to a product pure state $\varphi_{1}\otimes\varphi_{2}$ 
 where both $\varphi_{1}$ and $\varphi_{2}$ are pure. As $\varphi$ is not equivalent to $\varphi\circ\Theta$
 at least $\varphi_{1}$ or $\varphi_{2}$ is not equivalent to $\varphi_{1}\circ\Theta$ and $\varphi_{2}\circ\Theta$.
\par
If $\varphi_{1}$ and $\varphi_{1}\circ\Theta$  are equivalent,  there exists a $\Theta$ invariant pure state $\tilde{\varphi}_{1}$ 
equivalent to $\varphi_{1}$ and the argument of (i) applies to show our claim.
\par
 If $\varphi_{1}$ and $\varphi_{1}\circ\Theta$ are disjoint, due to Lemma \ref{lmm:30.9} ,
$\pi_{\varphi_{1}}(\mathfrak  A_{L}^{P (+)})^{\dprime} = \pi_{\varphi_{1}}(\mathfrak  A_{L}^{P})^{\dprime}$ is type I
and hence $\pi_{\psi}(\mathfrak  A_{L}^{(+)})^{\dprime} $ is type I and so is $\pi_{\psi}(\mathfrak  A_{L})^{\dprime} $.
\\
{\it End of Proof}
\bigskip
\par
All the staffs being ready, we summarize what we have shown in this section, which  automatically clarifies
 the relationship between string order and the index we defined.
 \par
 We started with a gapped ground state $\psi$ which is pure, and translationally invariant (and hence $\Theta$ invariant).
 Via Jordan-Wigner transformation, we obtain a  translationally invariant $\Theta$ invariant ground state $\psi^{P}$.
\par
If $ind_{\bfZ_{2}} \psi =1 $, $\psi^{P}$ is a gapped pure ground state ,  has no long range order.
For any $Q_{1} ,Q_{2} \in \mathfrak A^{P (-)}_{loc}$
$$\lim_{k\to \infty} \psi^{P} (Q_{1}\tau_{k}(Q_{2}) ) = 0$$
which implies
$$\lim_{k\to \infty} \psi^{P} (Q_{1} S [a , b+k] \tau_{k}(Q_{2}) ) = 0$$
for some $a$ and $b$   ( No String Order)
\par
If $ind_{\bfZ_{2}} \psi = -1 $, $\psi^{P}$ is still a ground state but not pure.
We have two possibilities 
\\
(a) $\varphi$ is a translationally invariant gapped pure ground state $\varphi$ of a Pauli spin system such that
$\varphi$ and $\varphi\circ\Theta$ are disjoint.
There exist  $Q_{1} ,Q_{2} \in \mathfrak A^{P (-)}_{loc}$ such that
$$ \lim_{k\to \infty} \psi^{P} (Q_{1}\tau_{k}(Q_{2}) ) \ne 0$$
Hence going back to $\mathfrak A$ string order exists
 $$\lim_{k\to \infty} \psi^{P} (Q_{1} S [a , b+k] \tau_{k}(Q_{2}) ) \ne 0$$
 (b) $\varphi$ is a periodic gapped ground state  of a Pauli spin system such that
  $$\varphi\circ\tau_{1} = \varphi\circ\Theta , \quad\varphi\circ\tau_{1} \ne \varphi .$$
There exist  $Q_{1} ,Q_{2} \in \mathfrak A^{P (-)}_{loc}$ such that
$$ \lim_{k\to \infty} \psi^{P} (Q_{1}\tau_{2k}(Q_{2}) ) \ne 0$$
String order exists
 $$\lim_{k\to \infty} \psi^{P} (Q_{1} S [a , b+ 2k] \tau_{k}(Q_{2}) ) \ne 0$$
Apparently the converse is correct. For example, if string order exists, $\psi^{P}$ is not pure and $ind_{\bfZ_{2}} \psi = -1 $. 
If no string order exists,  there is no long range order for $\psi^{P}$ and $ind_{\bfZ_{2}} \psi = 1 $.
This shows Theorem \ref{th:1.5}.
\bigskip
\\
{\bf Acknowledgment}
We would like to thank Y.Ogata (Univ. of Tokyo, Japan) for discussion. 
We would like to thank Hal Tasaki (Gakushin Univ. Japan) for showing us 
 the manuscript of his monograph prior to publication and for useful information.

\end{document}